\documentclass[reqno,11pt]{amsart}%
\usepackage{graphicx}
\usepackage{pstricks}
\usepackage{amsmath}
\usepackage{amsmath}
\usepackage{amsxtra}%
\usepackage{amsfonts}%
\usepackage{amssymb}
\theoremstyle{plain}
\newtheorem{theorem}{Theorem}[section]
\newtheorem{lemma}[theorem]{Lemma}
\newtheorem{proposition}[theorem]{Proposition}
\newtheorem{corollary}[theorem]{Corollary}
\newtheorem{conjecture}[theorem]{Conjecture}
\newtheorem{problem}[theorem]{Problem}
\theoremstyle{definition}
\newtheorem{definition}[theorem]{Definition}

\theoremstyle{remark}
\newtheorem{remark}[theorem]{Remark}

\numberwithin{equation}{section}
\begin{document}
\title{$k$-hyponormality of multivariable weighted shifts}
\author{Ra\'{u}l E. Curto}
\address{Department of Mathematics, The University of Iowa, Iowa City, Iowa 52242}
\email{rcurto@math.uiowa.edu}
\urladdr{http://www.math.uiowa.edu/\symbol{126}rcurto/}
\author{Sang Hoon Lee}
\address{Department of Mathematics, The University of Iowa, Iowa City, Iowa 52242}
\email{shlee@math.skku.ac.kr}
\urladdr{}
\author{Jasang Yoon}
\address{Department of Mathematics, Iowa State University, Ames, Iowa 50011}
\email{jyoon@iastate.edu}
\urladdr{http://www.public.iastate.edu/\symbol{126}jyoon/}
\thanks{The first named author was partially supported by NSF Grants DMS-0099357 and
DMS-0400741. The second named author was supported by the Post-doctoral
Fellowship Program of KOSEF }

\begin{abstract}
We characterize joint $k$-hyponormality for $2$-variable weighted shifts.
\ Using this characterization we construct a family of examples which
establishes and illustrates the gap between $k$-hyponormality and
$(k+1)$-hyponormality for each $k\geq1$. \ As a consequence, we obtain an
abstract solution to the Lifting Problem for Commuting Subnormals.

\end{abstract}
\maketitle

\section{\label{Int}Notation and Preliminaries}

The Lifting Problem for Commuting Subnormals asks for necessary and sufficient
conditions for a pair of subnormal operators on Hilbert space to admit
commuting normal extensions. \ It is well known that the commutativity of the
pair is necessary but not sufficient (\cite{Abr}, \cite{Lu1}, \cite{Lu2},
\cite{Lu3}), and it has recently been shown that the joint hyponormality of
the pair is necessary but not sufficient \cite{CuYo1}. \ In this paper we
provide an abstract answer to the Lifting Problem, by stating and proving a
multivariable analogue of the Bram-Halmos criterion for subnormality, and then
showing concretely that no matter how $k$-hyponormal a pair might be, it may
still fail to be subnormal. \ To do this, we obtain a matricial
characterization of $k$-hyponormality for multivariable weighted shifts, which
extends that found in \cite{bridge} for joint hyponormality.

Let $\mathcal{H}$ be a complex Hilbert space and let $\mathcal{B}%
(\mathcal{H})$ denote the algebra of bounded linear operators on $\mathcal{H}%
$. $\ $For $S,T\in\mathcal{B}(\mathcal{H})$ let $[S,T]:=ST-TS$. We say that an
$n$ -tuple $\mathbf{T}=(T_{1},\cdots,T_{n})$ of operators on $\mathcal{H}$ is
(jointly) \textit{hyponormal} if the operator matrix
\begin{equation}
\lbrack\mathbf{T}^{\ast},\mathbf{T]:=}\left(
\begin{array}
[c]{llll}%
\lbrack T_{1}^{\ast},T_{1}] & [T_{2}^{\ast},T_{1}] & \cdots & [T_{n}^{\ast
},T_{1}]\\
\lbrack T_{1}^{\ast},T_{2}] & [T_{2}^{\ast},T_{2}] & \cdots & [T_{n}^{\ast
},T_{2}]\\
\text{ \thinspace\thinspace\quad}\vdots & \text{ \thinspace\thinspace\quad
}\vdots & \ddots & \text{ \thinspace\thinspace\quad}\vdots\\
\lbrack T_{1}^{\ast},T_{n}] & [T_{2}^{\ast},T_{n}] & \cdots & [T_{n}^{\ast
},T_{n}]
\end{array}
\right)  \label{hypodef}%
\end{equation}
is positive on the direct sum of $n$ copies of $\mathcal{H}$ (cf. \cite{Ath} ,
\cite{CMX}). The $n$-tuple $\mathbf{T}$ is said to be \textit{normal} if
$\mathbf{T}$ is commuting and each $T_{i}$ is normal, and $\mathbf{T}$ is
\textit{subnormal }if $\mathbf{T}$ is the restriction of a normal $n$-tuple to
a common invariant subspace. Clearly, normal $\Rightarrow$ subnormal
$\Rightarrow$ hyponormal. \ Moreover, the restriction of a hyponormal
$n$-tuple to an invariant subspace is again hyponormal. \ The Bram-Halmos
criterion states that an operator $T\in\mathcal{B}(\mathcal{H})$ is subnormal
if and only if the $k$ -tuple $(T,T^{2},\cdots,T^{k})$ is hyponormal for all
$k\geq1$. \ 

For $\alpha\equiv\{\alpha_{n}\}_{n=0}^{\infty}$ a bounded sequence of positive
real numbers (called \textit{weights}), let $W_{\alpha}:\ell^{2}%
(\mathbb{Z}_{+})\rightarrow\ell^{2}(\mathbb{Z}_{+})$ be the associated
unilateral weighted shift, defined by $W_{\alpha}e_{n}:=\alpha_{n}e_{n+1}%
\;($all $n\geq0)$, where $\{e_{n}\}_{n=0}^{\infty}$ is the canonical
orthonormal basis in $\ell^{2}(\mathbb{Z}_{+}).$ The moments of $\alpha$ are
given as
\[
\gamma_{k}\equiv\gamma_{k}(\alpha):=\left\{
\begin{array}
[c]{cc}%
1 & \text{if }k=0\\
\alpha_{0}^{2}\cdots\alpha_{k-1}^{2} & \text{if }k>0
\end{array}
\right\}  .
\]
It is easy to see that $W_{\alpha}$ is never normal, and that it is hyponormal
if and only if $\alpha_{0}\leq\alpha_{1}\leq\cdots$. Similarly, consider
double-indexed positive bounded sequences $\alpha_{\mathbf{k}},\beta
_{\mathbf{k}}\in\ell^{\infty}(\mathbb{Z}_{+}^{2})$, $\mathbf{k}\equiv
(k_{1},k_{2})\in\mathbb{Z}_{+}^{2}:=\mathbb{Z}_{+}\times\mathbb{Z}_{+}$ and
let $\ell^{2}(\mathbb{Z}_{+}^{2})$ be the Hilbert space of square-summable
complex sequences indexed by $\mathbb{Z}_{+}^{2}$. \ We define the
$2$-\textit{variable weighted shift} $\mathbf{T}\equiv(T_{1},T_{2})$\ by
\[
T_{1}e_{\mathbf{k}}:=\alpha_{\mathbf{k}}e_{\mathbf{k+}\varepsilon_{1}}%
\]%
\[
T_{2}e_{\mathbf{k}}:=\beta_{\mathbf{k}}e_{\mathbf{k+}\varepsilon_{2}},
\]
where $\mathbf{\varepsilon}_{1}:=(1,0)$ and $\mathbf{\varepsilon}_{2}:=(0,1)$.
Clearly,
\begin{equation}
T_{1}T_{2}=T_{2}T_{1}\Longleftrightarrow\beta_{\mathbf{k+}\varepsilon_{1}%
}\alpha_{\mathbf{k}}=\alpha_{\mathbf{k+}\varepsilon_{2}}\beta_{\mathbf{k}%
}\;\;(\text{all }\mathbf{k}). \label{commuting}%
\end{equation}
In an entirely similar way one can define multivariable weighted shifts.
\ Trivially, a pair of unilateral weighted shifts $W_{\alpha}$ and $W_{\beta}$
gives rise to a $2$-variable weighted shift $\mathbf{T}\equiv(T_{1},T_{2})$,
if we let $\alpha_{(k_{1},k_{2})}:=\alpha_{k_{1}}$ and $\beta_{(k_{1},k_{2}%
)}:=\beta_{k_{2}}\;$(all $k_{1},k_{2}\in\mathbb{Z}_{+}$). \ In this case,
$\mathbf{T}$ is subnormal (resp. hyponormal) if and only if so are $T_{1}$ and
$T_{2}$; in fact, under the canonical identification of $\ell^{2}%
(\mathbb{Z}_{+}^{2})$ with $\ell^{2}(\mathbb{Z}_{+})\bigotimes\ell
^{2}(\mathbb{Z}_{+})$, we have $T_{1}\cong I\bigotimes W_{\alpha}$ and
$T_{2}\cong W_{\beta}\bigotimes I$, so $\mathbf{T}$ is also doubly commuting.
For this reason, we do not focus attention on shifts of this type, and use
them only when the above mentioned triviality is desirable or needed.

We now recall a well known characterization of subnormality for single
variable weighted shifts, due to C. Berger (cf. \cite[III.8.16]{Con}):
$W_{\alpha}$ is subnormal if and only if there exists a probability measure
$\xi$ supported in $[0,\left\|  W_{\alpha}\right\|  ^{2}]$ such that
$\gamma_{k}(\alpha):=\alpha_{0}^{2}\cdots\alpha_{k-1}^{2}=\int t^{k}%
\;d\xi(t)\;\;(k\geq1)$. \ If $W_{\alpha}$ is subnormal, and if for $h\geq1$ we
let $\mathcal{M}_{h}:=\bigvee\{e_{n}:n\geq h\}$ denote the invariant subspace
obtained by removing the first $h$ vectors in the canonical orthonormal basis
of $\ell^{2}(\mathbb{Z}_{+})$, then the Berger measure of $W_{\alpha
}|_{\mathcal{M}_{h}}$ is $\frac{1}{\gamma_{h}} t^{h}d\xi(t)$.

We also recall the notion of moment of order $\mathbf{k}$ for a pair
$(\alpha,\beta)$ satisfying (\ref{commuting}). \ Given $\mathbf{k}%
\in\mathbb{Z}_{+}^{2}$, the moment of $(\alpha,\beta)$ of order $\mathbf{k}$
is
\[
\gamma_{\mathbf{k}}\equiv\gamma_{\mathbf{k}}(\alpha,\beta):=\left\{
\begin{array}
[c]{cc}%
1 & \text{if }\mathbf{k}=0\\
\alpha_{(0,0)}^{2}\cdot...\cdot\alpha_{(k_{1}-1,0)}^{2} & \text{if }k_{1}%
\geq1\text{ and }k_{2}=0\\
\beta_{(0,0)}^{2}\cdot...\cdot\beta_{(0,k_{2}-1)}^{2} & \text{if }%
k_{1}=0\text{ and }k_{2}\geq1\\
\alpha_{(0,0)}^{2}\cdot...\cdot\alpha_{(k_{1}-1,0)}^{2}\cdot\beta_{(k_{1}%
,0)}^{2}\cdot...\cdot\beta_{(k_{1},k_{2}-1)}^{2} & \text{if }k_{1}\geq1\text{
and }k_{2}\geq1
\end{array}
\right\}  .
\]
We remark that, due to the commutativity condition (\ref{commuting}),
$\gamma_{\mathbf{k}}$ can be computed using any nondecreasing path from
$(0,0)$ to $(k_{1},k_{2})$. \ Moreover, $\mathbf{T}$ is subnormal if and only
if there is a regular Borel probability measure $\mu$ defined on the
$2$-dimensional rectangle $R=[0,a_{1}]\times\lbrack0,a_{2}]$ ($a_{i}:=\left\|
T_{i}\right\|  ^{2}$) such that $\gamma_{\mathbf{k}}=\iint_{R}\mathbf{t}%
^{\mathbf{k}}d\mu(\mathbf{t}):=\iint_{R}t_{1}^{k_{1}}t_{2}^{k_{2}}\;d\mu
(t_{1},t_{2})$ \ $($all $\mathbf{k\in}\mathbb{Z}_{+}^{2}$) \cite{JeLu}.

\textit{Acknowledgment}. \ The authors are deeply indebted to the referee for
several helpful suggestions. \ Many of the examples in this paper were
obtained using calculations with the software tool \textit{Mathematica
\cite{Wol}.}

\section{Main Results}

We recall some useful notation. \ For $n\geq0,$ let $m:=\frac{(n+1)(n+2)}{2}.$
\ For $A\in M_{m}(\mathbb{R}),$ we denote the successive rows and columns
according to the following lexicographic ordering: $1,x,y,x^{2},yx,y^{2},$
$\cdots,x^{n},yx^{n-1},\cdots,y^{n-1}x,y^{n}$ \cite{tcmp6}. \ For $0\leq
i+j\leq n,0\leq l+k\leq n,$ we denote the entry of $A\in M_{m}(\mathbb{R})$ in
row $y^{j}x^{i}$ and column $y^{k}x^{l}$ by $A_{(i,j)(l,k)}.$ \ In the
notation $0\leq i+j\leq n$ it will always be understood that $i,j\geq0.$ \ For
$0\leq i+j\leq n,0\leq l+k\leq n,$ $(a_{(i,j)(l,k)})_{_{0\leq l+k\leq
n}^{0\leq i+j\leq n}}$ denotes an $m\times m$ matrix and $(a_{(i,j)(l,k)}%
)_{_{1\leq l+k\leq n}^{1\leq i+j\leq n}}$ denotes the associated
$(m-1)\times(m-1)$ matrix obtained by deleting the first row and column.

For a subnormal $2$-variable weighted shift $\mathbf{T}\equiv(T_{1},T_{2})$,
it is clear that each component $T_{i}$ must be subnormal. \ For instance,
$T_{1}\cong\bigoplus_{j=0}^{\infty}W_{\alpha^{(j)}}$, where $\alpha_{i}%
^{(j)}:=\alpha_{(i,j)}$, so that $W_{\alpha^{(j)}}$ has associated Berger
measure $d\nu_{j}(t_{1}):=\frac{1}{\gamma_{(0,j)}}\int_{[0,a_{2}]}t_{2}%
^{j}d\Phi_{t_{1}}(t_{2})$, where $d\mu(t_{1},t_{2})\equiv d\Phi_{t_{1}}%
(t_{2})d\eta(t_{1})$ is the canonical disintegration of the Berger measure
$\mu$ by vertical slices \cite{CuYo2}. \ On the other hand, if we only know
that each of $T_{1}$, $T_{2}$ is subnormal, and that they commute, the
following problem is natural.

\begin{problem}
\label{problem} (Lifting Problem for Commuting Subnormals) \ Find necessary
and sufficient conditions on $T_{1}$ and $T_{2}$ to guarantee the subnormality
of $\mathbf{T}\equiv(T_{1},T_{2})$.
\end{problem}

It is well known that the above mentioned necessary conditions do not suffice
(cf.\cite{bridge}). \ In terms of the \textit{marginal} measures for $T_{1}$
and $T_{2}$, the problem can be phrased as a reconstruction-of-measure
problem, that is, under what conditions on the single variable measures
$\{\nu_{j}\}_{j=0}^{\infty}$ and $\{\omega_{i}\}_{i=0}^{\infty}$ associated
with $T_{1}$ and $T_{2}$, respectively, does there exist a $2$-variable
measure $\mu$ correctly interpolating all the powers $t_{1}^{k_{1}}%
t_{2}^{k_{2}}\;(k_{1},k_{2}\geq0)$ ? \ We also recall that a pair
$\mathbf{S}=(S_{1},S_{2})$ of commuting subnormal operators is called
\textit{polynomially subnormal} if $p(\mathbf{S})$ is subnormal for all
$2$-variable polynomials $p\in\mathbb{C}[z_{1},z_{2}]$. \ In \cite{Fra}, it
was shown that a polynomial subnormal tuple is a subnormal tuple. \ Using this
fact, we can give an abstract answer to Problem \ref{problem}. \ First we need
a definition.

\begin{definition}
A commuting pair $\mathbf{T}\equiv(T_{1},T_{2})$ is called $k$%
\textit{-hyponormal} if $\mathbf{T}(k):=(T_{1},T_{2},T_{1}^{2},T_{2}T_{1},$
$T_{2}^{2}$, $\cdots,T_{1}^{k},T_{2}T_{1}^{k-1},\cdots,T_{2}^{k})$ is
hyponormal, or equivalently
\[
([(T_{2}^{q}T_{1}^{p})^{\ast},T_{2}^{n}T_{1}^{m}])_{_{1\leq p+q\leq k}^{1\leq
m+n\leq k}}\geq0.
\]
\end{definition}

Clearly, subnormal $\Rightarrow$ $(k+1)$-hyponormal $\Rightarrow$
$k$-hyponormal for every $k\geq1,$ and of course $1$-hyponormality agrees with
the usual definition of joint hyponormality. \ 

We now present our multivariable version of the Bram-Halmos criterion for
subnormality. \ When combined with Theorem \ref{examplethm} below, Theorem
\ref{equivthm} provides an abstract answer to Problem \ref{problem}, by
showing that no matter how $k$-hyponormal the pair $\mathbf{T}$ might be, it
may still fail to be subnormal.

\begin{theorem}
\label{equivthm}Let $\mathbf{T}\equiv(T_{1},T_{2})$ be a commuting pair of
operators on a Hilbert space $\mathcal{H}$. \ The following statements are
equivalent.\newline (i) $\mathbf{T}$ is subnormal.\newline (ii) $\mathbf{T}%
(k)$ is subnormal for all $k\in\mathbb{Z}_{+}$.\newline (iii) $\mathbf{T}$ is
$k$-hyponormal for all $k\in\mathbb{Z}_{+}$.
\end{theorem}

In the single variable case, there are useful criteria for $k$-hyponormality
(\cite{QHWS}, \cite{CLL2}); for $2$-variable weighted shifts, a simple
criterion for joint hyponormality was given in(\cite{bridge}). \ We now
present a new characterization of $k$-hyponormality for $2$-variable weighted
shifts; this generalizes a result in (\cite{bridge}).

\begin{theorem}
\label{thm}Let $\mathbf{T}\equiv(T_{1},T_{2})$ be a $2$-variable weighted
shift with weight sequences $\alpha\equiv\{\alpha_{\mathbf{k}}\}$ and
$\beta\equiv\{\beta_{\mathbf{k}}\}$. \ The following statements are
equivalent.\newline (a) $\mathbf{T}$ is $k$-hyponormal.\newline (b)
$((T_{2}^{n}T_{1}^{m})^{\ast}[(T_{2}^{q}T_{1}^{p})^{\ast},T_{2}^{n}T_{1}%
^{m}](T_{2}^{q}T_{1}^{p}))_{_{1\leq p+q\leq k}^{1\leq m+n\leq k}}\geq
0$.\newline (c) $(\left\langle [(T_{2}^{q}T_{1}^{p})^{\ast},T_{2}^{n}T_{1}%
^{m}]e_{\mathbf{u}+(m,n)},e_{\mathbf{u}+(p,q)}\right\rangle )_{_{1\leq p+q\leq
k}^{1\leq m+n\leq k}}\geq0$ for all $\mathbf{u}\in\mathbb{Z}_{+}^{2}%
$.\newline (d) $(\gamma_{\mathbf{u}}\gamma_{\mathbf{u}+(m,n)+(p,q)}%
-\gamma_{\mathbf{u}+(m,n)}\gamma_{\mathbf{u}+(p,q)})_{_{1\leq p+q\leq
k}^{1\leq m+n\leq k}}\geq0$ for all $\mathbf{u}\in\mathbb{Z}_{+}^{2}%
$.\newline (e) $M_{\mathbf{u}}(k):=(\gamma_{\mathbf{u}+(m,n)+(p,q)})_{_{0\leq
p+q\leq k}^{0\leq m+n\leq k}}\geq0$ for all $\mathbf{u}\in\mathbb{Z}_{+}^{2}$.
\ (For a subnormal pair $\mathbf{T}$, the matrix $M_{\mathbf{u}}(k)$ is the
truncation of the moment matrix associated to the Berger measure of
$\mathbf{T}$.)
\end{theorem}

As an application of Theorem \ref{thm}, we build in Section \ref{appl} a
two-parameter family of $2$-variable weighted shifts (see Figure \ref{Figure
2} below), and we identify the precise parameter ranges that separate
hyponormality from $2$-hyponormality, $2$-hyponormality from $3$%
-hyponormality, etc., and $k$-hyponormality from subnormality. \ We believe
these are the first examples in the literature of commuting pairs of subnormal
operators which are $k$-hyponormal but not $(k+1)$-hyponormal. \ We record
this in the following result. \ First, we need some notation. \ For $0<y\leq
1$, let $x\equiv\{x_{n}\}_{n=0}^{\infty}$ where
\[
x_{n}:=\left\{
\begin{tabular}
[c]{ll}%
$y\sqrt{\frac{3}{4}},$ & $\text{if }n=0$\\
$\frac{\sqrt{(n+1)(n+3)}}{(n+2)},$ & $\text{if }n\geq1.$%
\end{tabular}
\ \right.
\]
(We shall see later (Proposition \ref{ex for k-hypo}) that $W_{x}\equiv
shift(x_{0},x_{1},\cdots)$ is subnormal.)

\begin{theorem}
\label{examplethm}For $0<a\leq\frac{1}{\sqrt{2}}$, the $2$-variable weighted
shift $\mathbf{T}$ given by Figure \ref{Figure 2} is\newline (i) hyponormal
$\Leftrightarrow0<y\leq\frac{\sqrt{32-48a^{4}}}{\sqrt{59-72a^{2}}};$%
\newline (ii) $k$-hyponormal $\Leftrightarrow$ $0<y\leq\sqrt
{\frac{\frac{(k+1)^{2}}{2k(k+2)}-a^{2}}{a^{4}-\frac{5}{2}a^{2}+\frac{(k+1)^{2}%
}{2k(k+2)}+\frac{2k^{2}+4k+3}{4(k+1)^{2}}}}\;\;(k\geq2);$\newline (iii)
subnormal $\Leftrightarrow0<y\leq\sqrt{\frac{1}{2-a^{2}}}$.\newline In
particular, $\mathbf{T}$ is hyponormal and not subnormal if and only if
$\sqrt{\frac{1}{2-a^{2}}}<y\leq\frac{\sqrt{32-48a^{4}}}{\sqrt{59-72a^{2}}}$.
\end{theorem}

\setlength{\unitlength}{1mm}  \psset{unit=1mm}  \begin{figure}[th]
\begin{center}
\begin{picture}(140,138)

\psline{->}%
(20,20)(135,20) \psline(20,40)(125,40)
\psline(20,60)(125,60) \psline(20,80)(125,80)
\psline(20,100)(125,100) \psline(20,120)(125,120)
\psline{->}%
(20,20)(20,135) \psline(40,20)(40,125)
\psline(60,20)(60,125) \psline(80,20)(80,125)
\psline(100,20)(100,125) \psline(120,20)(120,125)

\put(11,16){\footnotesize{$(0,0)$}%
}%

\put(35,16){\footnotesize{$(1,0)$}%
}%

\put(55,16){\footnotesize{$(2,0)$}%
}%

\put(78,16){\footnotesize{$\cdots$}%
}%

\put(95,16){\footnotesize{$(n,0)$}%
}%

\put(115,16){\footnotesize{$(n+1,0)$}%
}%

\put(27,21){\footnotesize{$x_0$}%
}%

\put(47,21){\footnotesize{$x_1$}%
}%

\put(67,21){\footnotesize{$x_2$}%
}%

\put(87,21){\footnotesize{$\cdots$}%
}%

\put(107,21){\footnotesize{$x_n$}%
}%

\put(124,21){\footnotesize{$\cdots$}%
}%

\put(27,41){\footnotesize{$a$}%
}%

\put(47,41){\footnotesize{$1$}%
}%

\put(67,41){\footnotesize{$1$}%
}%

\put(87,41){\footnotesize{$\cdots$}%
}%

\put(107,41){\footnotesize{$1$}%
}%

\put(124,41){\footnotesize{$\cdots$}%
}%

\put(27,61){\footnotesize{$a$}%
}%

\put(47,61){\footnotesize{$1$}%
}%

\put(67,61){\footnotesize{$1$}%
}%

\put(87,61){\footnotesize{$\cdots$}%
}%

\put(107,61){\footnotesize{$1$}%
}%

\put(124,61){\footnotesize{$\cdots$}%
}%

\put(27,81){\footnotesize{$\cdots$}%
}%

\put(47,81){\footnotesize{$\cdots$}%
}%

\put(67,81){\footnotesize{$\cdots$}%
}%

\put(87,81){\footnotesize{$\cdots$}%
}%

\put(107,81){\footnotesize{$\cdots$}%
}%

\put(124,81){\footnotesize{$\cdots$}%
}%

\put(27,101){\footnotesize{$a$}%
}%

\put(47,101){\footnotesize{$1$}%
}%

\put(67,101){\footnotesize{$1$}%
}%

\put(87,101){\footnotesize{$\cdots$}%
}%

\put(107,101){\footnotesize{$1$}%
}%

\put(124,101){\footnotesize{$\cdots$}%
}%

\put(27,121){\footnotesize{$a$}%
}%

\put(47,121){\footnotesize{$1$}%
}%

\put(67,121){\footnotesize{$1$}%
}%

\put(87,121){\footnotesize{$\cdots$}%
}%

\put(107,121){\footnotesize{$1$}%
}%

\put(124,121){\footnotesize{$\cdots$}%
}%

\psline{->}%
(70,10)(90,10) \put(79,6){$\rm{T}%
_1$}%

\put(11,38){\footnotesize{$(0,1)$}%
}%

\put(11,58){\footnotesize{$(0,2)$}%
}%

\put(14,78){\footnotesize{$\vdots$}%
}%

\put(11,98){\footnotesize{$(0,n)$}%
}%

\put(4,118){\footnotesize{$(0,n+1)$}%
}

\psline{->}%
(10, 70)(10,90) \put(5,80){$\rm{T}%
_2$}%

\put(20,28){\footnotesize{$y$}%
}%

\put(20,48){\footnotesize{$1$}%
}%

\put(20,68){\footnotesize{$1$}%
}%

\put(22,88){\footnotesize{$\vdots$}%
}%

\put(20,108){\footnotesize{$1$}%
}%

\put(22,128){\footnotesize{$\vdots$}%
}%

\put(40,28){\footnotesize{$\frac{ay}%
{x_0}%
$}%
}%

\put(40,48){\footnotesize{$1$}%
}%

\put(40,68){\footnotesize{$1$}%
}%

\put(42,88){\footnotesize{$\vdots$}%
}%

\put(40,108){\footnotesize{$1$}%
}%

\put(42,128){\footnotesize{$\vdots$}%
}%

\put(60,28){\footnotesize{$\frac{ay}%
{x_0x_1}%
$}%
}%

\put(60,48){\footnotesize{$1$}%
}%

\put(60,68){\footnotesize{$1$}%
}%

\put(62,88){\footnotesize{$\vdots$}%
}%

\put(60,108){\footnotesize{$1$}%
}%

\put(62,128){\footnotesize{$\vdots$}%
}%

\put(80,28){\footnotesize{$\frac{ay}%
{x_0x_1x_2}%
$}%
}%

\put(80,48){\footnotesize{$1$}%
}%

\put(80,68){\footnotesize{$1$}%
}%

\put(82,88){\footnotesize{$\vdots$}%
}%

\put(80,108){\footnotesize{$1$}%
}%

\put(82,128){\footnotesize{$\vdots$}%
}%

\put(100,28){\footnotesize{$\frac{ay}%
{x_0 \cdots x_{n-1}%
}%
$}%
}%

\put(100,48){\footnotesize{$1$}%
}%

\put(100,68){\footnotesize{$1$}%
}%

\put(102,88){\footnotesize{$\vdots$}%
}%

\put(100,108){\footnotesize{$1$}%
}%

\put(102,128){\footnotesize{$\vdots$}%
}%

\put(122,28){\footnotesize{$\vdots$}%
}%

\put(122,48){\footnotesize{$\vdots$}%
}%

\put(122,68){\footnotesize{$\vdots$}%
}%

\put(122,88){\footnotesize{$\vdots$}%
}%

\put(122,108){\footnotesize{$\vdots$}%
}%

\put(122,128){\footnotesize{$\vdots$}%
}%

\end{picture}
\end{center}
\caption{Weight diagram of the $2$-variable weighted shift in Theorem
\ref{examplethm} }%
\label{Figure 2}%
\end{figure}
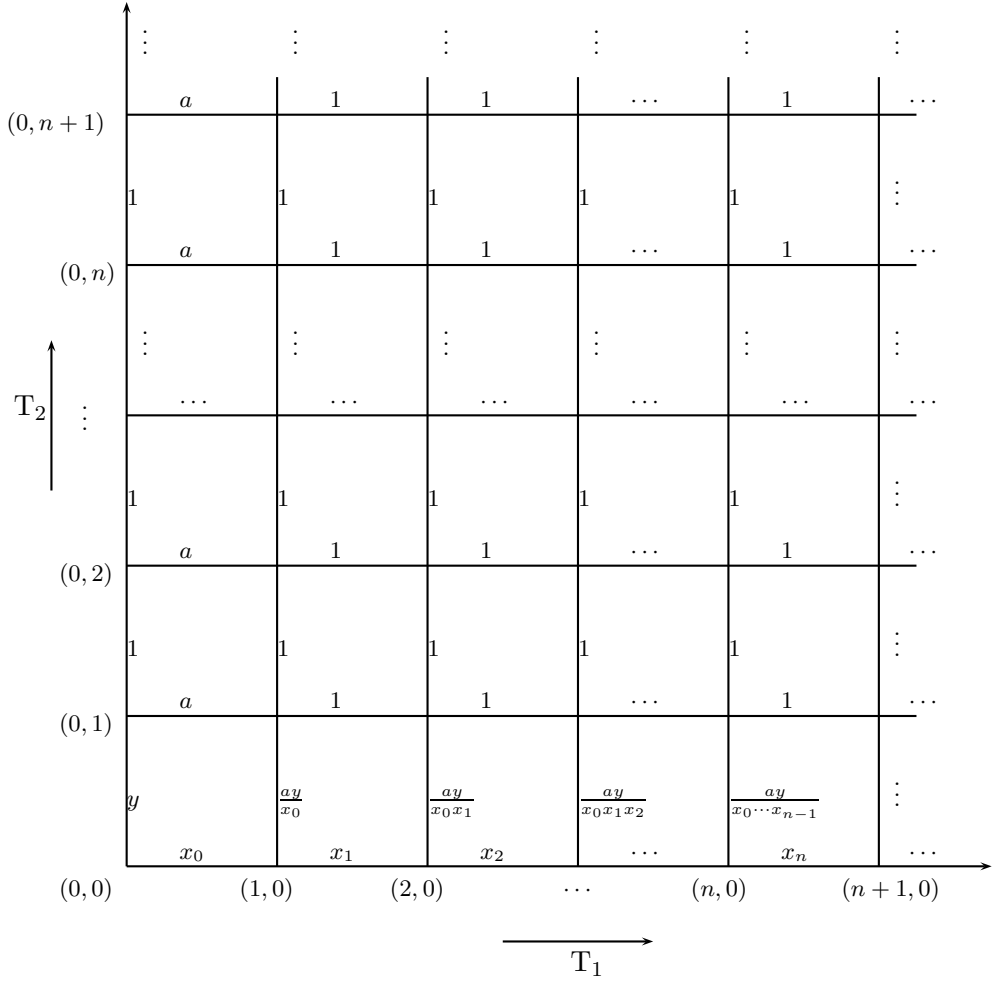

\begin{remark}
(i) \ Even for $1$-variable weighted shifts, it is generally difficult to
provide concrete parameterizations that separate $k$-hyponormality from
$(k+1)$-hyponormality (cf. \cite[Example 8]{CLL2}). \ That we can accomplish
the same separation for $2$-variable weighted shifts is an indication that the
condition in Theorem \ref{thm}(e) is sharp.\newline (ii) \ In \cite{CMX}, the
authors conjectured that if $\mathbf{T\equiv(}T_{1},T_{2})$ is a pair of
commuting subnormal operators, then $\mathbf{T}$ is subnormal if and only if
$\mathbf{T}$ is hyponormal. \ In \cite{CuYo1}, three different families of
examples were given of such pairs $\mathbf{T}$ for which hyponormality does
not imply subnormality. \ Thus, any of those examples can be used to disprove
the conjecture in \cite{CMX}. \ Theorem \ref{examplethm} gives a new family of
examples, with explicit parameter values to distinguish between $k$%
-hyponormality and $(k+1)$-hyponormality, and a fortiori between hyponormality
and subnormality.
\end{remark}

\section{Proofs of Theorems \ref{equivthm} and \ref{thm}}

\begin{proof}
[Proof of Theorem \ref{equivthm}]$(i)\Rightarrow(ii)$: \ Suppose
$\mathbf{T}\equiv(T_{1},T_{2})$ is subnormal, that is, $\mathbf{T}$ admits a
normal extension $\mathbf{N}\equiv(N_{1},N_{2})$ acting on a Hilbert space
$\mathcal{K}\supseteq\mathcal{H}$. \ The tuple $\mathbf{N}(k):=(N_{1}%
,N_{2},N_{1}^{2},N_{2}N_{1},N_{2}^{2}$, $\cdots,$ $N_{1}^{k},N_{2}N_{1}%
^{k-1},\cdots,N_{2}^{k})$ is also normal, so its restriction to $\mathcal{H}$,
$\mathbf{T}(k)$, is subnormal.

$(ii)\Rightarrow(iii)$: \ This is trivial.

$(iii)\Rightarrow(i)$: \ Suppose $\mathbf{T}(k)$ is hyponormal for all
$k\in\mathbb{Z}_{+}$, and let $p\in\mathbb{C}[z_{1},z_{2}]$. \ It follows that
$p(T_{1},T_{2})$ (as a single operator on $\mathcal{H}$) is $k$-hyponormal for
every $k\geq1$. \ By the Bram-Halmos criterion for single operators we then
see that $p(T_{1},T_{2})$ is subnormal. \ Finally, the main result in
\cite{Fra} implies that $\mathbf{T}$ is subnormal.
\end{proof}

We now give the proof of Theorem \ref{thm}, which we restate for the reader's convenience.

\begin{theorem}
Let $\mathbf{T}\equiv(T_{1},T_{2})$ be a $2$-variable weighted shift with
weight sequences $\alpha\equiv\{\alpha_{\mathbf{k}}\}$ and $\beta\equiv
\{\beta_{\mathbf{k}}\}$. \ The following statements are equivalent.\newline
(a) $\mathbf{T}$ is $k$-hyponormal. \newline (b) $((T_{2}^{n}T_{1}^{m})^{\ast
}[(T_{2}^{q}T_{1}^{p})^{\ast},T_{2}^{n}T_{1}^{m}](T_{2}^{q}T_{1}%
^{p}))_{_{1\leq p+q\leq k}^{1\leq m+n\leq k}}\geq0$. \newline (c)
$(\left\langle [(T_{2}^{q}T_{1}^{p})^{\ast},T_{2}^{n}T_{1}^{m}]e_{\mathbf{u}%
+(m,n)},e_{\mathbf{u}+(p,q)}\right\rangle )_{_{1\leq p+q\leq k}^{1\leq m+n\leq
k}}\geq0$ for all $\mathbf{u}\in\mathbb{Z}_{+}^{2}$. \newline (d)
$(\gamma_{\mathbf{u}}\gamma_{\mathbf{u}+(m,n)+(p,q)}-\gamma_{\mathbf{u}%
+(m,n)}\gamma_{\mathbf{u}+(p,q)})_{_{1\leq p+q\leq k}^{1\leq m+n\leq k}}\geq0$
for all $\mathbf{u}\in\mathbb{Z}_{+}^{2}$. \newline (e) $M_{\mathbf{u}%
}(k):=(\gamma_{\mathbf{u}+(m,n)+(p,q)})_{_{0\leq p+q\leq k}^{0\leq m+n\leq k}%
}\geq0$ for all $\mathbf{u}\in\mathbb{Z}_{+}^{2}$. \ 
\end{theorem}

\begin{proof}
$(a)\Leftrightarrow(b):$ \ Let
\[
A:=([(T_{2}^{q}T_{1}^{p})^{\ast},T_{2}^{n}T_{1}^{m}])_{_{1\leq p+q\leq
k}^{1\leq m+n\leq k}}%
\]
and
\[
B:=((T_{2}^{n}T_{1}^{m})^{\ast}[(T_{2}^{q}T_{1}^{p})^{\ast},T_{2}^{n}T_{1}%
^{m}](T_{2}^{q}T_{1}^{p}))_{_{1\leq p+q\leq k}^{1\leq m+n\leq k}}.
\]
Note that
\[
B=\text{$\left(
\begin{array}
[c]{cccccc}%
T_{1}^{\ast} &  &  &  &  & \\
& T_{2}^{\ast} &  &  &  & \\
&  & \ddots &  &  & \\
&  &  & (T_{1}^{k})^{\ast} &  & \\
&  &  &  & \ddots & \\
&  &  &  &  & (T_{2}^{k})^{\ast}%
\end{array}
\right)  $}A\left(
\begin{array}
[c]{cccccc}%
T_{1} &  &  &  &  & \\
& T_{2} &  &  &  & \\
&  & \ddots &  &  & \\
&  &  & T_{1}^{k} &  & \\
&  &  &  & \ddots & \\
&  &  &  &  & T_{2}^{k}%
\end{array}
\right)  .
\]
Thus it is easily seen that $(a)\Longrightarrow(b).$ For the converse
implication, let
\[
\mathcal{M}:=T_{1}(\mathcal{H})\oplus T_{2}(\mathcal{H})\oplus T_{1}%
^{2}(\mathcal{H})\oplus T_{2}T_{1}(\mathcal{H})\oplus T_{2}^{2}(\mathcal{H}%
)\oplus\cdots\oplus T_{1}^{k}(\mathcal{H})\oplus\cdots\oplus T_{2}%
^{k}(\mathcal{H})
\]
where $\mathcal{H}:=\ell^{2}(\mathbb{Z}_{+}^{2}).$ For
\begin{align*}
x  &  :=T_{1}(h_{(1,0)})\oplus T_{2}(h_{(0,1)})\oplus T_{1}^{2}(h_{(2,0)}%
)\oplus T_{2}T_{1}(h_{(1,1)})\oplus T_{2}^{2}(h_{(0,2)})\oplus\cdots\\
&  \oplus T_{1}^{k}(h_{(k,0)})\oplus\cdots\oplus T_{2}^{k}(h_{(0,k)})
\end{align*}
in $\mathcal{M}$ (here $h:=h_{(1,0)}\oplus h_{(0,1)}\oplus\cdots\oplus
h_{(0,k)}$), we have
\[
\left\langle Ax,x\right\rangle =\sum_{_{1\leq p+q\leq k}^{1\leq m+n\leq k}%
}\left\langle [(T_{2}^{q}T_{1}^{p})^{\ast},T_{2}^{n}T_{1}^{m}](T_{2}^{q}%
T_{1}^{p})h_{(p,q)},T_{2}^{n}T_{1}^{m}h_{(m,n)}\right\rangle =\left\langle
Bh,h\right\rangle .
\]
Thus, $\left\langle Bh,h\right\rangle \geq0$ implies that the compression of
$A$ to $\mathcal{M}$ is positive. \ Since $T_{1}$ and $T_{2}$ are weighted
shifts, we have $[(T_{2}^{q}T_{1}^{p})^{\ast},T_{2}^{n}T_{1}^{m}](T_{2}%
^{q}T_{1}^{p})(\mathcal{H})$ $\subseteq T_{2}^{n}T_{1}^{m}(\mathcal{H})$, and
it follows that $\mathcal{M}$ is invariant for $A$; since $A$ is selfadjoint,
it then follows that $\mathcal{M}$ reduces $A$. \ Therefore, the proof will be
completed once we show that the compression of $A$ to ${\mathcal{M}}^{\perp}$
is also positive. \ Let $\{e_{(p,q)}\}_{p,q\in\mathbb{Z}_{+}}$ be the
canonical orthonormal basis of $\mathcal{H}$. \ Note that
\[
{\mathcal{M}}^{\perp}=L_{(1,0)}\oplus L_{(0,1)}\oplus\cdots\oplus
L_{(k,0)}\oplus\cdots\oplus L_{(0,k)},
\]
where $L_{(i,j)}$ is the span of $\{e_{(p,q)}\}_{_{q\in\mathbb{Z}_{+}}^{0\leq
p\leq(i-1)}}$ and $\{e_{(r,s)}\}_{_{0\leq s\leq(j-1)}^{r\in\mathbb{Z}_{+}}}$.
$\ $Let
\[
e_{(\mathbf{u},\mathbf{m})}:=0\oplus0\oplus\cdots\oplus0\oplus\overset
{y^{j}x^{i}}{e_{\mathbf{u}}}\oplus0\oplus\cdots\oplus0\in\mathcal{H}^{m}%
\]
where $\mathbf{m}=(i,j).$ \ (The notation $0\oplus0\oplus\cdots\oplus
0\oplus\overset{y^{j}x^{i}}{\mathbf{v}}\oplus0\oplus\cdots\oplus0$ indicates
that the vector $\mathbf{v}$ appears in the $y^{j}x^{i}$-th summand.)
\ Observe that $\{e_{(p,q)(r,s)}\}_{_{q\in\mathbb{Z}_{+}}^{_{1\leq p\leq
(r-1)}^{1\leq r+s\leq k}}}\cup\{e_{(p,q)(r,s)}\}_{_{p\in\mathbb{Z}_{+}%
}^{_{1\leq p\leq(s-1)}^{1\leq r+s\leq k}}}$ is an orthonormal basis for
${\mathcal{M}}^{\perp}$. $\ $Now, for $1\leq i+j\leq k$, we define the
subspace $\mathcal{K}_{(i,j)}$ of ${\mathcal{M}}^{\perp}$ as follows:
\[
\mathcal{K}{_{(i,j)}}:=\left\{
\begin{tabular}
[c]{l}%
$\underset{n}{\mathbf{\vee}}\{e_{((0,n),(i,j))},e_{((1,n),(i+1,j))}%
,\cdots,e_{((k-(i+j),n),(k-j,j))}\},$ if $i>j$\\
$\underset{n}{\vee}\{e_{((n,0),(i,j))},e_{((n,1),(i,j+1))},\cdots
,e_{((n,k-(i+j)),(i,k-i))}\},$ if $i<j$\\
$\underset{n}{\vee}\{e_{((0,n),(i,i))},e_{((n,0),(i,i))},\cdots
,e_{((k-2i,n),(k-i,i))},e_{((n,k-2i),(i,k-i))}\},$ if $i=j.$%
\end{tabular}
\ \ \right.
\]
It is easily seen that $\mathcal{K}_{(i,j)}\perp\mathcal{K}_{(s,t)}$ for
$(i,j)\neq(s,t)$ and ${\mathcal{M}}^{\perp}=\mathcal{K}_{(1,0)}\oplus
\mathcal{K}_{(0,1)}\oplus\ldots\oplus\mathcal{K}_{(k,0)}\oplus\mathcal{K}%
_{(0,k)}.$ \ Note that $\mathcal{K}_{(i,j)}$ is invariant for $A$ and hence
$\mathcal{K}_{(i,j)}$ reduces $A$. \ We must then show that the compression of
$A$ to $\mathcal{K}_{(i,j)}$ is also positive. \ Note that each vector $h$ in
$\mathcal{K}_{(i,j)}(i>j)$ has the form
\begin{align*}
&  0\oplus\cdots\oplus0\oplus\overset{y^{j}x^{i}}{c(i,j)e_{(0,\ell)}}%
\oplus0\oplus\cdots\oplus0\oplus\overset{y^{j}x^{i+1}}{c(i+1,j)T_{1}%
e_{(0,\ell)}}\oplus0\oplus\cdots\oplus0\\
&  \oplus\overset{y^{j}x^{k-j}}{c(k-j,j)T_{1}^{k-(i+j)}e_{(0,\ell)}}%
\oplus0\oplus\cdots\oplus0
\end{align*}
for some scalars $c(i,j),c(i+1,j),\cdots,c(k-j,j)$ and $\ell\in\mathbb{Z}%
_{+}.$ \ Thus,%
\begin{equation}
\left\{
\begin{tabular}
[c]{l}%
$\left\langle Ah,h\right\rangle =\sum_{_{\substack{i\leq p\leq k-j,\ q=j}%
}^{^{\substack{i\leq m\leq k-j,\ n=j}}}}\left\langle [(T_{2}^{q}T_{1}%
^{p})^{\ast},T_{2}^{n}T_{1}^{m}]c(p,q)T_{1}^{p-i}e_{(0,\ell)},c(m,n)T_{1}%
^{m-i}e_{(0,\ell)}\right\rangle $\\
$=\sum_{_{i\leq p\leq k-j}^{^{i\leq m\leq k-j}}}c(p,j)\overline{c(m,j)}%
\left\langle T_{1}^{\ast^{(m-i)}}[(T_{2}^{j}T_{1}^{p})^{\ast},T_{2}^{j}%
T_{1}^{m}]T_{1}^{p-i}e_{(0,\ell)},e_{(0,\ell)}\right\rangle $\\
$=\sum_{_{i\leq p\leq k-j}^{^{i\leq m\leq k-j}}}c(p,j)\overline{c(m,j)}%
\left\langle (T_{2}^{j}T_{1}^{p+m-i})^{\ast}T_{2}^{j}T_{1}^{p+m-j}e_{(0,\ell
)},e_{(0,\ell)}\right\rangle $\\
$=\left\langle M\left(
\begin{array}
[c]{c}%
c(i,j)e_{(0,\ell)}\\
c(i+1,j)e_{(0,\ell)}\\
\vdots\\
c(k-j,j)e_{(0,\ell)}%
\end{array}
\right)  ,\left(
\begin{array}
[c]{c}%
c(i,j)e_{(0,\ell)}\\
c(i+1,j)e_{(0,\ell)}\\
\vdots\\
c(k-j,j)e_{(0,\ell)}%
\end{array}
\right)  \right\rangle ,$%
\end{tabular}
\right\}  \ \ \ \label{neweq}%
\end{equation}
$\ $

where%
\[
M:=\left(
\begin{array}
[c]{cccc}%
(T_{2}^{j}T_{1}^{i})^{\ast}(T_{2}^{j}T_{1}^{i}) & (T_{2}^{j}T_{1}^{i+1}%
)^{\ast}(T_{2}^{j}T_{1}^{i+1}) & \cdots & (T_{2}^{j}T_{1}^{k-j})^{\ast}%
(T_{2}^{j}T_{1}^{k-j})\\
(T_{2}^{j}T_{1}^{i+1})^{\ast}(T_{2}^{j}T_{1}^{i+1}) & (T_{2}^{j}T_{1}%
^{i+2})^{\ast}(T_{2}^{j}T_{1}^{i+2}) & \cdots & (T_{2}^{j}T_{1}^{k-j+1}%
)^{\ast}(T_{2}^{j}T_{1}^{k-j+1})\\
\vdots & \vdots & \ddots & \vdots\\
(T_{2}^{j}T_{1}^{k-j})^{\ast}(T_{2}^{j}T_{1}^{k-j}) & \cdots & \cdots &
(T_{2}^{j}T_{1}^{2(k-j)-i})^{\ast}(T_{2}^{j}T_{1}^{2(k-j)-i})
\end{array}
\right)  .
\]
In (\ref{neweq}) above, the third equality (from line 2 to line 3) follows
from the fact that
\[
\left\langle T_{2}^{j}T_{1}^{m}(T_{2}^{j}T_{1}^{p})^{\ast}T_{1}^{p-i}%
e_{(0,\ell)},e_{(0,\ell)}\right\rangle =0.
\]
Now,
\[
M=\left(
\begin{array}
[c]{ccc}%
(T_{2}^{j}T_{1}^{i})^{\ast} & 0 & 0\\
0 & \ddots & 0\\
0 & 0 & (T_{2}^{j}T_{1}^{i})^{\ast}%
\end{array}
\right)  M^{\prime}\left(
\begin{array}
[c]{ccc}%
T_{2}^{j}T_{1}^{i} & 0 & 0\\
0 & \ddots & 0\\
0 & 0 & T_{2}^{j}T_{1}^{i}%
\end{array}
\right)
\]
where
\[
M^{\prime}:=\left(
\begin{array}
[c]{cccc}%
I & T_{1}^{\ast}T_{1} & \cdots & (T_{1}^{k-(i+j)})^{\ast}T_{1}^{k-(i+j)}\\
T_{1}^{\ast}T_{1} & (T_{1}^{2})^{\ast}T_{1}^{2} & \cdots & (T_{1}%
^{k-(i+j)+1})^{\ast}T_{1}^{k-(i+j)+1}\\
\vdots & \vdots & \ddots & \vdots\\
(T_{1}^{k-(i+j)})^{\ast}T_{1}^{k-(i+j)} & (T_{1}^{k-(i+j)+1})^{\ast}%
T_{1}^{k-(i+j)+1} & \cdots & (T_{1}^{2(k-(i+j))})^{\ast}T_{1}^{2(k-(i+j))}%
\end{array}
\right)  .
\]
Now, by Smul'jan's Theorem \cite{Smu}, $M^{\prime}\equiv(M_{uv}^{\prime
})_{u,v=0}^{k-(i+j)}\geq0$ if and only if $Q\equiv(Q_{uv})_{u,v=1}%
^{k-(i+j)}\geq0$, where
\[
Q_{uv}:=(T_{1}^{u+v})^{\ast}T_{1}^{u+v}-(T_{1}^{u})^{\ast}T_{1}^{u}(T_{1}%
^{v})^{\ast}T_{1}^{v}.
\]
Now, observe that $Q_{uv}=(T_{1}^{u})^{\ast}[(T_{1}^{v})^{\ast},T_{1}%
^{u}]T_{1}^{v}$, so that
\[
Q=\text{$\left(
\begin{array}
[c]{ccc}%
T_{1}^{\ast} &  & \\
& \ddots & \\
&  & (T_{1}^{k})^{\ast}%
\end{array}
\right)  \left(
\begin{array}
[c]{ccc}%
\lbrack T_{1}^{\ast},T_{1}] & \cdots & [(T_{1}^{k-(i+j)})^{\ast},T_{1}]\\
\vdots & \ddots & \vdots\\
\lbrack(T_{1})^{\ast},T_{1}^{k-(i+j)}] & \cdots & [(T_{1}^{k-(i+j)})^{\ast
},T_{1}^{k-(i+j)}]
\end{array}
\right)  \left(
\begin{array}
[c]{ccc}%
T_{1} &  & \\
& \ddots & \\
&  & T_{1}^{k}%
\end{array}
\right)  $}%
\]
It is now easy to see that $Q$ is a submatrix of $B$. \ Thus, if $B\geq0$,
then $M^{\prime}\geq0$ and hence $\left\langle Ah,h\right\rangle \geq0$ for
all $h\in\mathcal{K}_{(i,j)}$ with $i>j$. \ On the other hand, if $i<j$, then
each vector $h$ in $\mathcal{K}_{(i,j)}$ has the form
\begin{align*}
&  0\oplus\cdots\oplus0\oplus\overset{y^{j}x^{i}}{c(i,j)e_{(\ell,0)}}%
\oplus0\oplus\cdots\oplus0\oplus\overset{y^{j+1}x^{i}}{c(i,j+1)T_{2}%
e_{(\ell,0)}}\oplus0\oplus\cdots\\
&  \oplus0\oplus\overset{y^{k-i}x^{i}}{c(i,k-i)T_{2}^{k-(i+j)}e_{(\ell,0)}%
}\oplus0\oplus\cdots\oplus0,
\end{align*}
for some scalars $c(i,j),c(i,j+1),\cdots,c(i,k-i)$ and $\ell\in\mathbb{Z}%
_{+}.$ \ An analogous argument shows that $\left\langle Ah,h\right\rangle
\geq0$ for all $h\in\mathcal{K}_{(i,j)}$ with $i<j$. \ Finally, if $i=j$, then
each vector $h$ in $\mathcal{K}_{(i,i)}$ has the form
\begin{align*}
&  0\oplus\cdots\oplus0\oplus\overset{y^{i}x^{i}}{c(i,i)e_{(s,t)}}%
\oplus0\oplus\cdots\oplus0\oplus\overset{y^{i}x^{i+1}}{c(i+1,i)T_{1}e_{(s,t)}%
}\\
&  \oplus0\oplus\cdots\oplus0\oplus\overset{y^{i+1}x^{i}}{c(i,i+1)T_{2}%
e_{(s,t)}}\oplus0\oplus\cdots\oplus0\oplus\overset{y^{i}x^{k-i}}%
{c(k-i,i)T_{1}^{k-2i}e_{(s,t)}}\\
&  \oplus0\oplus\cdots\oplus0\oplus\overset{y^{k-i}x^{i}}{c(i,k-i)T_{2}%
^{k-2i}e_{(s,t)}}\oplus0\oplus\cdots\oplus0
\end{align*}
for some scalars $c(i,i),c(i+1,i),c(i,i+1)\cdots,c(i,k-i)$ and $s,t\in
\mathbb{Z}_{+}$ with $st=0.$ \ Define
\[
\sigma_{(p,q)}:=\left\{
\begin{array}
[c]{cc}%
1 & \text{if }p\geq q\\
2 & \text{if }p<q,
\end{array}
\right.
\]
and let $M(p,q):=\max\{p,q\}$. \ We then have%
\[%
\begin{tabular}
[c]{l}%
$\left\langle Ah,h\right\rangle =\underset{_{_{\substack{2i\leq p+q\leq
k,\ p=i\ \text{or}\ q=i}}^{^{\substack{2i\leq m+n\leq k,\ n=i\ \text{or}%
\ m=i}}}}}{\sum c(p,q)\overline{c(m,n)}}\left\langle [(T_{2}^{q}T_{1}%
^{p})^{\ast},T_{2}^{n}T_{1}^{m}]T_{\sigma_{(p,q)}}^{M(p,q)-i}e_{(s,t)}%
,T_{\sigma_{(m,n)}}^{M(m,n)-i}e_{(s,t)}\right\rangle $\\
$=\underset{_{_{\substack{2i\leq p+q\leq k,\ p=i\ \text{or}\ q=i}%
}^{^{\substack{2i\leq m+n\leq k,\ n=i\ \text{or}\ m=i}}}}}{\sum
c(p,q)\overline{c(m,n)}}\left\langle (T_{\sigma_{(m,n)}}^{M(m,n)-i})^{\ast
}[(T_{2}^{q}T_{1}^{p})^{\ast},T_{2}^{n}T_{1}^{m}]T_{\sigma_{(p,q)}}%
^{M(p,q)-i}e_{(s,t)},e_{(s,t)}\right\rangle .$%
\end{tabular}
\ \ \ \
\]
Since $(T_{2}^{q}T_{1}^{p})^{\ast}T_{\sigma_{(p,q)}}^{M(p,q)-i}e_{(s,t)}=0,$
we have%
\[%
\begin{tabular}
[c]{l}%
$\left\langle Ah,h\right\rangle =\underset{_{_{\substack{2i\leq p+q\leq
k,\ p=i\ \text{or}\ q=i}}^{^{\substack{2i\leq m+n\leq k,\ n=i\ \text{or}%
\ m=i}}}}}{\sum c(p,q)\overline{c(m,n)}}\left\langle (T_{\sigma_{(m,n)}%
}^{M(m,n)-i})^{\ast}(T_{2}^{q}T_{1}^{p})^{\ast}T_{2}^{n}T_{1}^{m}%
T_{\sigma_{(p,q)}}^{M(p,q)-i}e_{(s,t)},e_{(s,t)}\right\rangle $\\
$=\left\langle M\left(
\begin{array}
[c]{c}%
c(i,i)e_{(s,t)}\\
c(i+1,i)e_{(s,t)}\\
\vdots\\
c(i,k-i)e_{(s,t)}%
\end{array}
\right)  ,\left(
\begin{array}
[c]{c}%
c(i,i)e_{(s,t)}\\
c(i+1,i)e_{(s,t)}\\
\vdots\\
c(i,k-i)e_{(s,t)}%
\end{array}
\right)  \right\rangle ,$%
\end{tabular}
\ \ \ \
\]
where
\[
M:=\left(
\begin{array}
[c]{cccc}%
(T_{2}^{i}T_{1}^{i})^{\ast}(T_{2}^{i}T_{1}^{i}) & (T_{2}^{i}T_{1}^{i+1}%
)^{\ast}(T_{2}^{i}T_{1}^{i+1}) & \cdots & (T_{2}^{k-i}T_{1}^{i})^{\ast}%
(T_{2}^{k-i}T_{1}^{i})\\
(T_{2}^{i}T_{1}^{i+1})^{\ast}(T_{2}^{i}T_{1}^{i+1}) & (T_{2}^{i}T_{1}%
^{i+2})^{\ast}(T_{2}^{i}T_{1}^{i+2}) & \cdots & (T_{2}^{k-i}T_{1}^{i+1}%
)^{\ast}(T_{2}^{k-i}T_{1}^{i+1})\\
\vdots & \vdots & \ddots & \vdots\\
(T_{2}^{k-i}T_{1}^{i})^{\ast}(T_{2}^{k-i}T_{1}^{i}) & \cdots & \cdots &
(T_{2}^{k-i+1}T_{1}^{i})^{\ast}(T_{2}^{k-i+1}T_{1}^{i})
\end{array}
\right)  .
\]
However,
\[
M=\left(
\begin{array}
[c]{ccc}%
(T_{2}^{i}T_{1}^{i})^{\ast} & 0 & 0\\
0 & \ddots & 0\\
0 & 0 & (T_{2}^{i}T_{1}^{i})^{\ast}%
\end{array}
\right)  M^{\prime}\left(
\begin{array}
[c]{ccc}%
T_{2}^{i}T_{1}^{i} & 0 & 0\\
0 & \ddots & 0\\
0 & 0 & T_{2}^{i}T_{1}^{i}%
\end{array}
\right)  ,
\]
where
\[
M^{\prime}:=\left(
\begin{array}
[c]{cccc}%
I & T_{1}^{\ast}T_{1} & \cdots & (T_{2}^{k-2i})^{\ast}T_{2}^{k-2i}\\
T_{1}^{\ast}T_{1} & (T_{1}^{2})^{\ast}T_{1}^{2} & \cdots & (T_{2}^{k-2i}%
T_{1})^{\ast}T_{2}^{k-2i}T_{1}\\
\vdots & \vdots & \ddots & \vdots\\
(T_{2}^{k-2i})^{\ast}T_{2}^{k-2i} & \cdots & \cdots & (T_{2}^{2(k-2i)})^{\ast
}T_{2}^{2(k-2i)}%
\end{array}
\right)  .
\]
If $B\geq0$, then $M^{\prime}\geq0$ and hence $\left\langle Ah,h\right\rangle
\geq0$ for all $h\in\mathcal{K}_{(i,i)}.$

$(b)\Leftrightarrow(c):$ \ Note that $(T_{2}^{n}T_{1}^{m})^{\ast}[(T_{2}%
^{q}T_{1}^{p})^{\ast},T_{2}^{n}T_{1}^{m}](T_{2}^{q}T_{1}^{p})$ is a diagonal
operator. Thus, $B\geq0$ if and only if $\langle Be_{\mathbf{u}}%
,e_{\mathbf{u}}\rangle\geq0$ for all $\mathbf{u}\in\mathbb{Z}_{+}^{2}$ if and
only if $(c)$ holds.

$(b)\Leftrightarrow(d):$ \ Since $(T_{2}^{m}T_{1}^{n})^{\ast}[(T_{2}^{q}%
T_{1}^{p})^{\ast},T_{2}^{m}T_{1}^{n}](T_{2}^{q}T_{1}^{p})$ is a diagonal
operator whose $\mathbf{u}$-th diagonal entry is
\[
\frac{\gamma_{\mathbf{u}+(p,q)+(m,n)}}{\gamma_{\mathbf{u}}}-\frac{\gamma
_{\mathbf{u}+(p,q)}\gamma_{\mathbf{u}+(m,n)}}{\gamma_{\mathbf{u}}^{2}},
\]
we easily see that $(b)\Leftrightarrow(d).$

$(d)\Leftrightarrow(e):$ \ This is a straightforward application of Choleski's
algorithm \cite{Atk}.
\end{proof}

\section{Applications\label{appl}}

Unlike the single variable case, in which there is a clear separation between
hyponormality and subnormality (cf. \cite{RGWSII}, \cite{OTAMP},
\cite{CuLe3},\cite{CLL2}), much less is known about the multivariable case.
\ We will now construct an example which exhibits the gap between
$k$-hyponormality and $(k+1)$-hyponormality for each $k\geq1$, and gives
another counterexample to the following conjecture, recently answered in the
negative (\cite{CuYo1}).

\begin{conjecture}
\label{conjecture}(\cite{CMX}) Let $\mathbf{T\equiv(}T_{1},T_{2})$ be a pair
of commuting subnormal operators on $\mathcal{H}$. Then $\mathbf{T}$ is
subnormal if and only if $\mathbf{T}$ is hyponormal.
\end{conjecture}

We begin with:

\begin{proposition}
\label{ex for k-hypo}For $0<y\leq1$, let $x\equiv\{x_{n}\}_{n=0}^{\infty}$
where
\[
x_{n}:=\left\{
\begin{tabular}
[c]{ll}%
$y\sqrt{\frac{3}{4}},$ & $\text{if }n=0$\\
$\frac{\sqrt{(n+1)(n+3)}}{(n+2)},$ & $\text{if }n\geq1.$%
\end{tabular}
\right.
\]
Then $W_{x}\equiv shift(x_{0},x_{1},\cdots)$ is subnormal.
\end{proposition}

\begin{proof}
We need to find a regular Borel probability measure $\mu_{x}$ such that
$\gamma_{n}=\int s^{n}d\mu_{x}(s)\;\;(n\geq0)$. \ On the interval $[0,1]$,
consider $d\mu_{x}:=(1-y^{2})d\delta_{0}(s)+\frac{y^{2}}{2}ds+\frac{y^{2}}%
{2}d\delta_{1}(s)$. \ Then $\gamma_{0}=1$ and for $n\geq1,$
\begin{align*}
\gamma_{n}  &  \equiv x_{0}^{2}x_{1}^{2}x_{2}^{2}\cdots x_{n-1}^{2}\\
&  =y^{2}\frac{3}{2^{2}}\cdot\frac{2\cdot4}{3^{2}}\cdot\frac{3\cdot5}{4^{2}%
}\cdot\cdots\cdot\frac{n(n+2)}{(n+1)^{2}}\\
&  =\frac{(n+2)y^{2}}{2(n+1)}=\frac{y^{2}}{2}\cdot\frac{1}{n+1}+\frac{y^{2}%
}{2}=\int s^{n}d\mu_{x}(s).
\end{align*}
It follows that $W_{x}$ is subnormal, with Berger measure $\mu_{x}$.
\end{proof}

\setlength{\unitlength}{1mm}  \psset{unit=1mm}  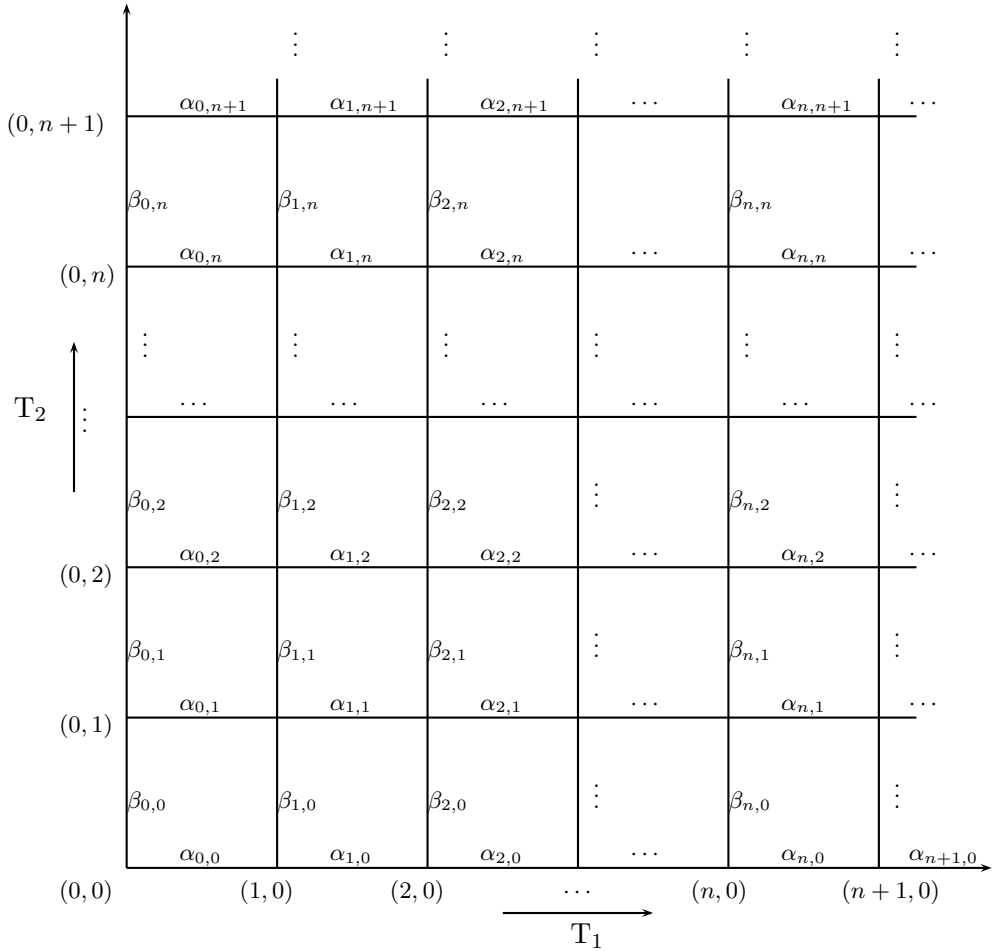
\begin{figure}[th]
\begin{center}
\begin{picture}(140,138)

\psline{->}%
(20,20)(135,20)
\psline(20,40)(125,40)
\psline(20,60)(125,60)
\psline(20,80)(125,80)
\psline(20,100)(125,100)
\psline(20,120)(125,120)
\psline{->}%
(20,20)(20,135)
\psline(40,20)(40,125)
\psline(60,20)(60,125)
\psline(80,20)(80,125)
\psline(100,20)(100,125)
\psline(120,20)(120,125)

\put(11,16){\footnotesize{$(0,0)$}%
}
\put(35,16){\footnotesize{$(1,0)$}%
}%

\put(55,16){\footnotesize{$(2,0)$}%
}%

\put(78,16){\footnotesize{$\cdots$}%
}%

\put(95,16){\footnotesize{$(n,0)$}%
}%

\put(115,16){\footnotesize{$(n+1,0)$}%
}%

\put(27,21){\footnotesize{$\alpha_{0,0}%
$}%
}%

\put(47,21){\footnotesize{$\alpha_{1,0}%
$}%
}%

\put(67,21){\footnotesize{$\alpha_{2,0}%
$}%
}%

\put(87,21){\footnotesize{$\cdots$}%
}%

\put(107,21){\footnotesize{$\alpha_{n,0}%
$}%
}%

\put(124,21){\footnotesize{$\alpha_{n+1,0}%
$}%
}%

\put(27,41){\footnotesize{$\alpha_{0,1}%
$}%
}%

\put(47,41){\footnotesize{$\alpha_{1,1}%
$}%
}%

\put(67,41){\footnotesize{$\alpha_{2,1}%
$}%
}%

\put(87,41){\footnotesize{$\cdots$}%
}%

\put(107,41){\footnotesize{$\alpha_{n,1}%
$}%
}%

\put(124,41){\footnotesize{$\cdots$}%
}%

\put(27,61){\footnotesize{$\alpha_{0,2}%
$}%
}%

\put(47,61){\footnotesize{$\alpha_{1,2}%
$}%
}%

\put(67,61){\footnotesize{$\alpha_{2,2}%
$}%
}%

\put(87,61){\footnotesize{$\cdots$}%
}%

\put(107,61){\footnotesize{$\alpha_{n,2}%
$}%
}%

\put(124,61){\footnotesize{$\cdots$}%
}%

\put(27,81){\footnotesize{$\cdots$}%
}%

\put(47,81){\footnotesize{$\cdots$}%
}%

\put(67,81){\footnotesize{$\cdots$}%
}%

\put(87,81){\footnotesize{$\cdots$}%
}%

\put(107,81){\footnotesize{$\cdots$}%
}%

\put(124,81){\footnotesize{$\cdots$}%
}%

\put(27,101){\footnotesize{$\alpha_{0,n}%
$}%
}%

\put(47,101){\footnotesize{$\alpha_{1,n}%
$}%
}%

\put(67,101){\footnotesize{$\alpha_{2,n}%
$}%
}%

\put(87,101){\footnotesize{$\cdots$}%
}%

\put(107,101){\footnotesize{$\alpha_{n,n}%
$}%
}%

\put(124,101){\footnotesize{$\cdots$}%
}%

\put(27,121){\footnotesize{$\alpha_{0,n+1}%
$}%
}%

\put(47,121){\footnotesize{$\alpha_{1,n+1}%
$}%
}%

\put(67,121){\footnotesize{$\alpha_{2,n+1}%
$}%
}%

\put(87,121){\footnotesize{$\cdots$}%
}%

\put(107,121){\footnotesize{$\alpha_{n,n+1}%
$}%
}%

\put(124,121){\footnotesize{$\cdots$}%
}%

\psline{->}%
(70,14)(90,14)
\put(79,10){$\rm{T}%
_1$}%

\put(11,38){\footnotesize{$(0,1)$}%
}%

\put(11,58){\footnotesize{$(0,2)$}%
}%

\put(14,78){\footnotesize{$\vdots$}%
}%

\put(11,98){\footnotesize{$(0,n)$}%
}%

\put(4,118){\footnotesize{$(0,n+1)$}%
}

\psline{->}%
(13,70)(13,90)
\put(5,80){$\rm{T}%
_2$}%

\put(20,28){\footnotesize{$\beta_{0,0}%
$}%
}%

\put(20,48){\footnotesize{$\beta_{0,1}%
$}%
}%

\put(20,68){\footnotesize{$\beta_{0,2}%
$}%
}%

\put(22,88){\footnotesize{$\vdots$}%
}%

\put(20,108){\footnotesize{$\beta_{0,n}%
$}%
}%

\put(40,28){\footnotesize{$\beta_{1,0}%
$}%
}%

\put(40,48){\footnotesize{$\beta_{1,1}%
$}%
}%

\put(40,68){\footnotesize{$\beta_{1,2}%
$}%
}%

\put(42,88){\footnotesize{$\vdots$}%
}%

\put(40,108){\footnotesize{$\beta_{1,n}%
$}%
}%

\put(42,128){\footnotesize{$\vdots$}%
}%

\put(60,28){\footnotesize{$\beta_{2,0}%
$}%
}%

\put(60,48){\footnotesize{$\beta_{2,1}%
$}%
}%

\put(60,68){\footnotesize{$\beta_{2,2}%
$}%
}%

\put(62,88){\footnotesize{$\vdots$}%
}%

\put(60,108){\footnotesize{$\beta_{2,n}%
$}%
}%

\put(62,128){\footnotesize{$\vdots$}%
}%

\put(82,28){\footnotesize{$\vdots$}%
}%

\put(82,48){\footnotesize{$\vdots$}%
}%

\put(82,68){\footnotesize{$\vdots$}%
}%

\put(82,88){\footnotesize{$\vdots$}%
}%

\put(82,128){\footnotesize{$\vdots$}%
}%

\put(100,28){\footnotesize{$\beta_{n,0}%
$}%
}%

\put(100,48){\footnotesize{$\beta_{n,1}%
$}%
}%

\put(100,68){\footnotesize{$\beta_{n,2}%
$}%
}%

\put(102,88){\footnotesize{$\vdots$}%
}%

\put(100,108){\footnotesize{$\beta_{n,n}%
$}%
}%

\put(102,128){\footnotesize{$\vdots$}%
}%

\put(122,28){\footnotesize{$\vdots$}%
}%

\put(122,48){\footnotesize{$\vdots$}%
}%

\put(122,68){\footnotesize{$\vdots$}%
}%

\put(122,88){\footnotesize{$\vdots$}%
}%

\put(122,128){\footnotesize{$\vdots$}%
}%

\end{picture}
\end{center}
\caption{Weight diagram of the 2-variable weighted shift in Lemma
\ref{backext}}%
\label{Figure 1}%
\end{figure}

To recall the following result, we need some notation and terminology from
\cite{CuYo1}. \ Given a probability measure $\mu$ on $X\times Y\equiv
\mathbb{R}_{+}\times\mathbb{R}_{+}$, and assuming that $\frac{1}{t}\in
L^{1}(\mu)$, the \emph{extremal measure} $\mu_{ext}$ (which is also a
probability measure) on $\mathbb{R}_{+}\times\mathbb{R}_{+}$ is given by
$d\mu_{ext}(s,t):=(1-\delta_{0}(t))\frac{1}{t\left\|  \frac{1}{t}\right\|
_{L^{1}(\mu)}}d\mu(s,t)$. \ On the other hand, the \textit{marginal measure}
$\mu^{X}$ is given by $\mu^{X}:=\mu\circ\pi_{X}^{-1}$, where $\pi_{X}:X\times
Y\rightarrow X$ is the canonical projection onto $X$. \ Thus, $\mu^{X}%
(E)=\mu(E\times Y)$, for every $E\subseteq X$. \ We observe that if $\mu$ is a
probability measure, then so is $\mu^{X}$.

\begin{lemma}
\label{backext}\cite[Proposition 3.10]{CuYo1} (Subnormal backward extension of
a $2$-variable weighted shift) \ Consider the $2$-variable weighted shift
$\mathbf{T}$ whose weight sequence is given by Figure \ref{Figure 1}, and let
$\mathcal{M}$ be the subspace associated with indices $\mathbf{k}$ with
$k_{2}\geq1$. \ Assume that $\mathbf{T|}_{\mathcal{M}}$ is subnormal with
Berger measure $\mu_{\mathcal{M}}$ and that the weighted shift $W_{0}$ with
weight sequence $(\alpha_{00},\alpha_{10},\cdots)$ is subnormal with Berger
measure $\nu$. \ Then $\mathbf{T}$ is subnormal if and only if \newline (i)
$\ \frac{1}{t}\in L^{1}(\mu_{\mathcal{M}})$;\newline (ii) $\ \beta_{00}%
^{2}\leq(\left\|  \frac{1}{t}\right\|  _{L^{1}(\mu_{\mathcal{M}})})^{-1}%
$;\newline (iii) $\ \beta_{00}^{2}\left\|  \frac{1}{t}\right\|  _{L^{1}%
(\mu_{\mathcal{M}})}(\mu_{\mathcal{M}})_{ext}^{X}\leq\nu$.\newline Moreover,
if $\beta_{00}^{2}\left\|  \frac{1}{t}\right\|  _{L^{1}(\mu_{\mathcal{M}}%
)}=1,$ then $(\mu_{\mathcal{M}})_{ext}^{X}=\nu$. \ In the case when
$\mathbf{T}$ is subnormal, the Berger measure $\mu$ of $\mathbf{T}$ is given
by
\[
d\mu(s,t)=\beta_{00}^{2}\left\|  \frac{1}{t}\right\|  _{L^{1}(\mu
_{\mathcal{M}})}d(\mu_{\mathcal{M}})_{ext}(s,t)+[d\nu(s)-\beta_{00}%
^{2}\left\|  \frac{1}{t}\right\|  _{L^{1}(\mu_{\mathcal{M}})}d(\mu
_{\mathcal{M}})_{ext}^{X}(s)]d\delta_{0}(t).
\]
\end{lemma}

We are now ready to present our example of a nonsubnormal, hyponormal
commuting pair of subnormal weighted shifts. \ At the same time we will
exhibit concretely the gap between $k$-hyponormality and $(k+1)$-hyponormality
for each $k\geq1$. For $0<a\leq\frac{1}{\sqrt{2}}$, consider the $2$-variable
weighted shift given by Figure \ref{Figure 2}, where $x\equiv\{x_{n}%
\}_{n=0}^{\infty}$ is as in Proposition \ref{ex for k-hypo}.

\begin{proposition}
\label{propsub}The $2$-variable weighted shift $\mathbf{T}$ given by Figure
\ref{Figure 2} is subnormal if and only if $0<y\leq\sqrt{\frac{1}{2-a^{2}}}.$
\end{proposition}

\begin{proof}
Let $\mathcal{M}$ be the subspace of $\ell^{2}(\mathbb{Z}_{+}^{2})$ spanned by
the canonical orthonormal basis of $\ell^{2}(\mathbb{Z}_{+}^{2})$ except for
$e_{(0,0)},e_{(1,0)},\cdots,e_{(n,0)},\cdots$. Then from Figure \ref{Figure
2}, it is obvious that $\mathbf{T|}_{\mathcal{M}}\cong(I\otimes S_{a}%
,U_{+}\otimes I)$. \ (Recall that $S_{a}$ is the subnormal weighted shift
which has weight sequence $(a,1,1,\cdots)$ and Berger measure $(1-a^{2}%
)\delta_{0}+a^{2}\delta_{1}$, and $U_{+}\equiv S_{1}$ is the (unweighted)
unilateral shift.) \ Thus, $\mathbf{T|}_{\mathcal{M}}$ is subnormal with
Berger measure
\[
\mu_{\mathcal{M}}:=[(1-a^{2})\delta_{0}+a^{2}\delta_{1}]\times\delta_{1}.
\]
By Lemma \ref{backext},
\begin{align*}
\mathbf{T}\text{ is subnormal } &  \Leftrightarrow y^{2}\left\|  \frac{1}%
{t}\right\|  _{L^{1}(\mu_{\mathcal{M}})}(\mu_{\mathcal{M}})_{ext}^{X}\leq
\mu_{x}\\
&  \Leftrightarrow y^{2}[(1-a^{2})\delta_{0}+a^{2}\delta_{1}]\leq
(1-y^{2})\delta_{0}+\frac{y^{2}}{2}\lambda+\frac{y^{2}}{2}\delta_{1}\\
&  \text{(here }\lambda\text{ denotes Lebesgue measure on }[0,1]\text{)}\\
&  \Leftrightarrow y\leq\sqrt{\frac{1}{2-a^{2}}}\text{ and }a\leq
\frac{1}{\sqrt{2}}\\
&  \Leftrightarrow y\leq\sqrt{\frac{1}{2-a^{2}}}\;\;\text{(recall that }%
a\leq\frac{1}{\sqrt{2}}\text{ is being assumed).}%
\end{align*}
\end{proof}

\begin{proposition}
\label{prophyp}The $2$-variable weighted shift $\mathbf{T}$ given by Figure
\ref{Figure 2} is hyponormal if and only if $0<y\leq\frac{\sqrt{32-48a^{4}}%
}{\sqrt{59-72a^{2}}}.$
\end{proposition}

\begin{proof}
By Theorem \ref{thm}(d), to show the joint hyponormality of $\mathbf{T}$ it is
enough to check that
\[
H_{\mathbf{k}}:=\left(
\begin{array}
[c]{cc}%
\gamma_{\mathbf{k}}\gamma_{\mathbf{k}+(2,0)}-\gamma_{\mathbf{k}+(1,0)}^{2} &
\gamma_{\mathbf{k}}\gamma_{\mathbf{k}+(1,1)}-\gamma_{\mathbf{k}+(1,0)}%
\gamma_{\mathbf{k}+(0,1)}\\
\gamma_{\mathbf{k}}\gamma_{\mathbf{k}+(1,1)}-\gamma_{\mathbf{k}+(1,0)}%
\gamma_{\mathbf{k}+(0,1)} & \gamma_{\mathbf{k}}\gamma_{\mathbf{k}%
+(0,2)}-\gamma_{\mathbf{k}+(0,1)}^{2}%
\end{array}
\right)  \geq0
\]
for all $\mathbf{k}\in\mathbb{Z}_{+}^{2}.$ Since $\mathbf{T|}_{\mathcal{M}}$
is subnormal (as noted in Proposition \ref{propsub}), it is also hyponormal,
so it remains to show that $H_{\mathbf{k}}\geq0$ for $\mathbf{k}%
=(0,0),(1,0),(2,0),\cdots,(n,0),\cdots$. A straightforward calculation shows
that
\begin{align*}
H_{(n,0)}  &  =\left(
\begin{array}
[c]{cccc}%
\gamma_{(n,0)}\gamma_{(n+2,0)}-\gamma_{(n+1,0)}^{2} & \gamma_{(n,0)}%
\gamma_{(n+1,1)}-\gamma_{(n+1,0)}\gamma_{(n,1)} &  & \\
\gamma_{(n,0)}\gamma_{(n+1,1)}-\gamma_{(n+1,0)}\gamma_{(n,1)} & \gamma
_{(n,0)}\gamma_{(n,2)}-\gamma_{(n,1)}^{2} &  &
\end{array}
\right) \\
&  =\gamma_{(n,0)}\left(
\begin{array}
[c]{cc}%
\gamma_{(n+2,0)}-{x_{n}^{2}}\gamma_{(n+1,0)} & \gamma_{(n+1,1)}-{x_{n}^{2}%
}\gamma_{(n,1)}\\
\gamma_{(n+1,1)}-{x_{n}^{2}}{\gamma_{(n,1)}} & \gamma_{(n,2)}-\frac{\gamma
_{(n,1)}^{2}}{\gamma_{(n,0)}}%
\end{array}
\right) \\
&  =\gamma_{(n,0)}\left(
\begin{array}
[c]{cc}%
(x_{0}\cdots x_{n+1})^{2}-x_{n}^{2}(x_{0}\cdots x_{n})^{2} & a^{2}y^{2}%
-x_{n}^{2}a^{2}y^{2}\\
a^{2}y^{2}-x_{n}^{2}a^{2}y^{2} & a^{2}y^{2}-\frac{(a^{2}y^{2})^{2}}%
{(x_{0}\cdots x_{n-1})^{2}}%
\end{array}
\right) \\
&  =y^{2}\gamma_{(n,0)}\left(
\begin{array}
[c]{cc}%
\frac{n+4}{2(n+3)}-\frac{n+3}{2(n+2)}\cdot\frac{(n+1)(n+3)}{(n+2)^{2}} &
a^{2}(1-\frac{(n+1)(n+3)}{(n+2)^{2}})\\
a^{2}(1-\frac{(n+1)(n+3)}{(n+2)^{2}}) & a^{2}(1-\frac{(n+2)}{2(n+1)}a^{2})
\end{array}
\right) \\
&  =y^{2}\gamma_{(n,0)}\left(
\begin{array}
[c]{cc}%
\frac{2n+5}{2(n+2)^{3}(n+3)} & \frac{1}{(n+2)^{2}}a^{2}\\
\frac{1}{(n+2)^{2}}a^{2} & a^{2}(1-\frac{(n+2)}{2(n+1)}a^{2})
\end{array}
\right)  ,
\end{align*}
which is positive because $0<a\leq\frac{1}{\sqrt{2}}$ for all $n\neq0.$ \ For
$\mathbf{k}=(0,0)$, we have
\[
H_{(0,0)}=\left(
\begin{array}
[c]{cc}%
(x_{0}x_{1})^{2}-(x_{0}^{2})^{2} & a^{2}y^{2}-x_{0}^{2}y^{2}\\
a^{2}y^{2}-x_{0}^{2}y^{2} & y^{2}-y^{4}%
\end{array}
\right)  .
\]
Since $(x_{0}x_{1})^{2}-(x_{0}^{2})^{2}>0$ and
\[
\det H_{(0,0)}=y^{4}(\frac{2}{3}-\frac{9}{16}y^{2})(1-y^{2})-(a^{2}%
-\frac{3}{4}y^{2})^{2}=\frac{y^{4}}{48}\{(72a^{2}-59)y^{2}+32-48a^{4}\},
\]
we obtain the desired result.
\end{proof}

In the following theorem, we summarize the results in Propositions
\ref{propsub} and \ref{prophyp}, and provide a new family of examples to
settle Conjecture \ref{conjecture} in the negative.

\begin{theorem}
\label{secondthm}If $\sqrt{\frac{1}{2-a^{2}}}<y\leq\frac{\sqrt{32-48a^{4}}%
}{\sqrt{59-72a^{2}}}$, then the $2$-variable weighted shift $\mathbf{T}$ given
by Figure \ref{Figure 2} is hyponormal but not subnormal.
\end{theorem}

For $k\geq2$ we let
\[
H_{k}(y):=\frac{1}{2}\left(
\begin{array}
[c]{cccccccc}%
\frac{2}{y^{2}} & \frac{3}{2} & 2 & 2a^{2} & \frac{4}{3} & \frac{5}{4} &
\cdots & \frac{k+2}{k+1}\\
\frac{3}{2} & \frac{4}{3} & 2a^{2} & 2a^{2} & \frac{5}{4} & \frac{6}{5} &
\cdots & \frac{k+3}{k+2}\\
2 & 2a^{2} & 2 & 2a^{2} & 2a^{2} & 2a^{2} & \cdots & 2a^{2}\\
2a^{2} & 2a^{2} & 2a^{2} & 2a^{2} & 2a^{2} & 2a^{2} & \cdots & 2a^{2}\\
\frac{4}{3} & \frac{5}{4} & 2a^{2} & 2a^{2} & \frac{6}{5} & \frac{7}{6} &
\cdots & \frac{k+4}{k+3}\\
\frac{5}{4} & \frac{6}{5} & 2a^{2} & 2a^{2} & \frac{7}{6} & \frac{8}{7} &
\cdots & \frac{k+5}{k+4}\\
\vdots & \vdots & \vdots & \vdots & \vdots & \vdots & \ddots & \vdots\\
\frac{k+2}{k+1} & \frac{k+3}{k+2} & 2a^{2} & 2a^{2} & \frac{k+4}{k+3} &
\frac{k+5}{k+4} & \cdots & \frac{2k+2}{2k+1}%
\end{array}
\right)  _{(k+3)\times(k+3)}.
\]

We then have:

\begin{theorem}
\label{k-hypo}The $2$-variable weighted shift $\mathbf{T}$ given by Figure
\ref{Figure 2} is $k$-hyponormal ($k\geq2$) if and only if $0<y\leq
D(k):=\sqrt{\frac{\frac{(k+1)^{2}}{2k(k+2)}-a^{2}}{a^{4}-\frac{5}{2}%
a^{2}+\frac{(k+1)^{2}}{2k(k+2)}+\frac{2k^{2}+4k+3}{4(k+1)^{2}}}}$.
\end{theorem}

\begin{remark}
Since $D(k+1)<D(k)$ for every $k\geq2$, it follows that for $D(k+1)<y\leq
D(k)$, the associated $2$-variable weighted shift $\mathbf{T}$ is
$k$-hyponormal but not $(k+1)$-hyponormal.
\end{remark}

\begin{proof}
By Theorem \ref{thm}(e), $\mathbf{T}$ is $k$-hyponormal if and only if
\[
M_{\mathbf{k}}(k)=(\gamma_{\mathbf{k}+(m,n)+(p,q)})_{_{0\leq p+q\leq k}^{0\leq
n+m\leq k}}\geq0,
\]
for all $\mathbf{k}\in\mathbb{Z}_{+}^{2}.$ By \cite[Theorem 5.2 and Remark
5.3]{CuYo1} we need to verify that $M_{\mathbf{k}}(k)\geq0$ for $\mathbf{k}%
=(0,0).$ A direct computation shows that this is equivalent to $H_{k}(y)\geq0$
and, in turn, equivalent to $\det H_{k}(y)\geq0$ and the fact that
$\mathbf{T}$ is $(k-1)$-hyponormal. Now let
\[
A_{k}:=\frac{1}{2}\left(
\begin{array}
[c]{cccccccc}%
\frac{4}{3} & 2a^{2} & 2a^{2} & \frac{5}{4} & \frac{6}{5} & \cdots &
\frac{k+3}{k+2} & \\
2a^{2} & 2 & 2a^{2} & 2a^{2} & 2a^{2} & \cdots & 2a^{2} & \\
2a^{2} & 2a^{2} & 2a^{2} & 2a^{2} & 2a^{2} & \cdots & 2a^{2} & \\
\frac{5}{4} & 2a^{2} & 2a^{2} & \frac{6}{5} & \frac{7}{6} & \cdots &
\frac{k+4}{k+3} & \\
\frac{6}{5} & 2a^{2} & 2a^{2} & \frac{7}{6} & \frac{8}{7} & \cdots &
\frac{k+5}{k+4} & \\
\vdots & \vdots & \vdots & \vdots & \vdots & \ddots & \vdots & \\
\frac{k+3}{k+2} & 2a^{2} & 2a^{2} & \frac{k+4}{k+3} & \frac{k+5}{k+4} & \cdots
& \frac{2k+2}{2k+1} &
\end{array}
\right)  _{(k+2)\times(k+2)}.
\]
Note that we can easily calculate $\det A_{k}$ and $\det H_{k}(1)$. Indeed,
observe that
\begin{equation}
\left\{
\begin{array}
[c]{l}%
\det A_{k}=a_{k}\cdot{a^{2}}(1-a^{2})(\frac{(k+1)^{2}}{2k(k+2)}-a^{2})\\
\det H_{k}(1)=a_{k}\cdot{a^{2}}(a^{2}-1)\{(1-a^{2})(\frac{1}{2}-a^{2}%
)+\frac{1}{4(k+1)^{2}}\}\\
a_{k}:=\frac{(1!2!\cdots(k-1)!)^{2}2!3!\cdots(k+1)!}{2^{k-1}(k+2)!(k+3)!\cdots
(2k+1)!}k(k+2).
\end{array}
\right.  \label{eq41}%
\end{equation}
On the other hand, by the cofactor expansion along the first row or the first
column, we have
\[
\det H_{k}(y)=\frac{1}{y^{2}}\det A_{k}+\det H_{k}(1)-\det A_{k}%
=(\frac{1}{y^{2}}-1)\det A_{k}+\det H_{k}(1).
\]
Since $\det A_{k}\geq0$,
\[
\det H_{k}(y)\geq0\Leftrightarrow y^{2}\leq\frac{\det A_{k}}{\det A_{k}-\det
H_{k}(1)}.
\]

Thus,
\[
\det H_{k}(y)\geq0\Leftrightarrow y^{2}\leq\frac{\frac{(k+1)^{2}}%
{2k(k+2)}-a^{2}}{a^{4}-\frac{5}{2}a^{2}+\frac{(k+1)^{2}}{2k(k+2)}%
+\frac{2k^{2}+4k+3}{4(k+1)^{2}}}.\
\]
Therefore, we see that $\mathbf{T}$ is $k$-hyponormal ($k\geq2$) if and only
if $0<y\leq\sqrt{\frac{\frac{(k+1)^{2}}{2k(k+2)}-a^{2}}{a^{4}-\frac{5}{2}%
a^{2}+\frac{(k+1)^{2}}{2k(k+2)}+\frac{2k^{2}+4k+3}{4(k+1)^{2}}}}$, as desired.
\end{proof}

\begin{corollary}
\label{final}Let $\mathbf{T}$ be the $2$-variable weighted shift given by
Figure \ref{Figure 2}, and assume that $\mathbf{T}$ is $k$-hyponormal for
every $k\geq2$. \ Then $0<y\leq\sqrt{\frac{1}{2-a^{2}}}$. \ 
\end{corollary}

\begin{proof}
We know that $0<y\leq D(k)$ for every $k\geq2$, so
\[
y\leq\lim_{k\rightarrow\infty}D(k)=\sqrt{\frac{\frac{1}{2}-a^{2}}%
{a^{4}-\frac{5}{2}a^{2}+1}}=\sqrt{\frac{\frac{1}{2}-a^{2}}{(a^{2}-\frac{1}%
{2})(a^{2}-2)}}=\sqrt{\frac{1}{2-a^{2}}},
\]
as desired.
\end{proof}

\begin{remark}
The results in Proposition \ref{propsub} and Corollary \ref{final} illustrate
Theorem \ref{equivthm} (the multivariable Bram-Halmos criterion); that is, the
pair $\mathbf{T}$ is subnormal if and only if it is $k$-hyponormal for every
$k\geq1$.
\end{remark}


\begin{thebibliography}{CuFi1}                                                                                            %

\bibitem[Abr]{Abr}M.B. Abrahamse, Commuting subnormal operators.
\textit{Illinois J. Math}. 22 (1978), 171--176.

\bibitem[Ath]{Ath}A. Athavale, On joint hyponormality of operators,
\textit{\ Proc. Amer. Math. Soc}. 103(1988), 417-423.

\bibitem[Atk]{Atk}K. Atkinson, \textit{Introduction to Numerical Analysis},
Wiley and Sons, 2nd. Ed. \ 1989.

\bibitem[Con]{Con}J. Conway, \textit{The Theory of Subnormal Operators,}
Mathematical Surveys and Monographs, vol. 36, Amer. Math. Soc., Providence, 1991.

\bibitem[Cu1]{bridge}R. Curto, Joint hyponormality: A bridge between
hyponormality and subnormality, \textit{Proc. Symposia Pure Math}. 51(1990), 69-91.

\bibitem[Cu2]{QHWS}R. Curto, Quadratically hyponormal weighted shifts,
\textit{Integral Equations Operator Theory} 13(1990), 49-66.

\bibitem[Cu3]{OTAMP}R. Curto, An operator-theoretic approach to truncated
moment problems, in \textit{Linear Operators}, Banach Center Publ., vol. 38,
1997; pp. 75-104.

\bibitem[CuFi1]{RGWSII}R. Curto and L. Fialkow, Recursively generated weighted
shifts and the subnormal completion problem, II, \textit{Integral Equations
Operator Theory}, 18(1994), 369-426.

\bibitem[CuFi2]{tcmp6}R. Curto and L. Fialkow, Solution of the singular
quartic moment problem, \textit{J. Operator Theory} 48(2002), 315-354.

\bibitem[CuLe]{CuLe3}R. Curto and W.Y. Lee, Towards a model theory for $2$
-hyponormal operators, \textit{Integral Equations Operator Theory} 44(2002), 290-315.

\bibitem[CLL]{CLL2}R. Curto, S.H. Lee and W.Y. Lee, A new criterion for
$k$-hyponormality via weak subnormality, \textit{Proc. Amer. Math. Soc}.
133(2005), 1805-1816.

\bibitem[CMX]{CMX}R. Curto, P. Muhly and J. Xia, Hyponormal pairs of commuting
operators, \textit{Operator Theory: Adv. Appl.} 35(1988), 1-22.

\bibitem[CuYo1]{CuYo1}R. Curto and J. Yoon, Jointly hyponormal pairs of
subnormal operators need not be jointly subnormal, Trans. Amer. Math. Soc., to appear.

\bibitem[CuYo2]{CuYo2}R. Curto and J. Yoon, Disintegration-of-measure
techniques for multivariable weighted shifts, preprint 2004.

\bibitem[Fra]{Fra}E. Franks, Polynomially subnormal operator tuples,
\textit{J. Operator Theory} 31(1994), 219-228.

\bibitem[JeLu]{JeLu}N.P. Jewell and A.R. Lubin, Commuting weighted shifts and
analytic function theory in several variables, \textit{J. Operator Theory }
1(1979), 207-223.

\bibitem[Lu1]{Lu1}A. Lubin, Weighted shifts and products of subnormal
operators. \textit{Indiana Univ. Math. J}. 26 (1977), 839--845.

\bibitem[Lu2]{Lu2}A. Lubin, Extensions of commuting subnormal operators, in
\textit{Hilbert space operators} (Proc. Conf., Calif. State Univ., Long Beach,
Calif., 1977), pp. 115--120, \textit{Lecture Notes in Math}. 693, Springer,
Berlin, 1978.

\bibitem[Lu3]{Lu3}A. Lubin, A subnormal semigroup without normal extension.
\textit{Proc. Amer. Math. Soc}. 68 (1978), 176--178.

\bibitem[Smu]{Smu}Ju. L. Smul'jan, An operator Hellinger integral,
\textit{Mat. Sb}. (N.S.) \textbf{49} (1959), 381--430 (Russian).

\bibitem[Wol]{Wol}Wolfram Research, Inc. \textit{Mathematica}, Version 4.2,
\textit{Wolfram Research Inc.}, Champaign, IL, 2002.
\end{thebibliography}
\end{document}